\documentclass[12pt,leqno,draft]{article}
\usepackage{amsfonts}
\usepackage{amsmath, amsthm, amsfonts, amssymb, color}
\usepackage{mathrsfs}
\usepackage{color}
\usepackage{stmaryrd}
\makeatletter
\newcommand{\Spvek}[2][r]{%
  \gdef\@VORNE{1}
  \left(\hskip-\arraycolsep%
    \begin{array}{#1}\vekSp@lten{#2}\end{array}%
  \hskip-\arraycolsep\right)}

\def\vekSp@lten#1{\xvekSp@lten#1;vekL@stLine;}
\def\vekL@stLine{vekL@stLine}
\def\xvekSp@lten#1;{\def\temp{#1}%
  \ifx\temp\vekL@stLine
  \else
    \ifnum\@VORNE=1\gdef\@VORNE{0}
    \else\@arraycr\fi%
    #1%
    \expandafter\xvekSp@lten
  \fi}
\makeatother

\newtheorem{theorem}{Theorem}[section]
\newtheorem{lemma}[theorem]{Lemma}

\newtheorem{example}[theorem]{Example}
\newtheorem{remark}[theorem]{Remark}
\theoremstyle{definition}
\newtheorem{definition}{Definition}[section]
\newcommand{\scr}[1]{\mathscr #1}
\definecolor{wco}{rgb}{0.5,0.2,0.3}

\numberwithin{equation}{section} \theoremstyle{remark}

\newcommand{\ua}{\uparrow}

\title{{\bf    Harnack Inequality and  Applications for SDEs Driven by $G$-Brownian motion}
}
\author{
{\bf     Fenfen Yang  }\\
\footnotesize{  Center for Applied Mathematics, Tianjin University, Tianjin 300072, China}\\
\footnotesize{  yangfenfen@tju.edu.cn}}
\begin{document}
\allowdisplaybreaks
\def\R{\mathbb R}  \def\ff{\frac} \def\ss{\sqrt} \def\B{\mathbf
B}
\def\N{\mathbb N} \def\kk{\kappa} \def\m{{\bf m}}
\def\ee{\varepsilon}\def\ddd{D^*}
\def\dd{\delta} \def\DD{\Delta} \def\vv{\varepsilon} \def\rr{\rho}
\def\<{\langle} \def\>{\rangle} \def\GG{\Gamma} \def\gg{\gamma}
  \def\nn{\nabla} \def\pp{\partial} \def\E{\mathbb E}
\def\d{\text{\rm{d}}} \def\bb{\beta} \def\aa{\alpha} \def\D{\scr D}
  \def\si{\sigma} \def\ess{\text{\rm{ess}}}
\def\beg{\begin} \def\beq{\begin{equation}}  \def\F{\scr F}
\def\Ric{\text{\rm{Ric}}} \def\Hess{\text{\rm{Hess}}}
\def\e{\text{\rm{e}}} \def\ua{\underline a} \def\OO{\Omega}  \def\oo{\omega}
 \def\tt{\tilde} \def\Ric{\text{\rm{Ric}}}
\def\cut{\text{\rm{cut}}} \def\P{\mathbb P} \def\ifn{I_n(f^{\bigotimes n})}
\def\C{\scr C}   \def\G{\scr G}   \def\aaa{\mathbf{r}}     \def\r{r}
\def\gap{\text{\rm{gap}}} \def\prr{\pi_{{\bf m},\varrho}}  \def\r{\mathbf r}
\def\Z{\mathbb Z} \def\vrr{\varrho} 
\def\L{\scr L}\def\Tt{\tt} \def\TT{\tt}\def\II{\mathbb I}
\def\i{{\rm in}}\def\Sect{{\rm Sect}}  \def\H{\mathbb H}
\def\M{\scr M}\def\Q{\mathbb Q} \def\texto{\text{o}} \def\LL{\Lambda}
\def\Rank{{\rm Rank}} \def\B{\scr B} \def\i{{\rm i}} \def\HR{\hat{\R}^d}
\def\to{\rightarrow}\def\l{\ell}\def\iint{\int}
\def\EE{\scr E}\def\no{\nonumber}
\def\A{\scr A}\def\V{\mathbb V}\def\osc{{\rm osc}}
\def\BB{\scr B}\def\Ent{{\rm Ent}}\def\3{\triangle}\def\H{\scr H}
\def\U{\scr U}\def\8{\infty}\def\1{\lesssim}\def\HH{\mathrm{H}}
 \def\T{\scr T}
\maketitle

\begin{abstract}  We  establish Harnack inequality and shift  Harnack inequality for stochastic differential equation driven by $G$-Brownian motion. As applications, the uniqueness of invariant linear  expectations and estimates on the $\sup$-kernel  are  investigated, where the
 $\sup$-kernel  is introduced in this paper for the first time.
\end{abstract} \noindent
 \vskip 2cm

\section{Introduction}
Since Wang \cite{W} introduced dimensional-free Harnack inequality for diffusions on Riemannian manifold, his Harnack inequality  has been  extensively investigated.
 His type Harnack inequality  as a powerful tool in the study of  functional inequalities (see \cite{Aida, RW1, RW2, W5, W6}), heat kernel estimates (see \cite{Gong}), high order eigenvalues (see  \cite{ W7, Kawabi}), transportation cost inequalities (see \cite{Bobkov}), and short-time behavior of transition probabilities (see \cite{AK, AZ, Kawabi}). To establish Wang's Harnack inequality,
 Wang and co-authors introduced the coupling by change of measures, see Wang \cite{W1} and references within for details.

On the other hand,
for the potential applications in uncertainty problems, risk measures and the superhedging in finance, the theory of nonlinear expectation has been developed. Especially, Peng \cite{peng2, peng1} established the fundamental theory of $G$-expectation theory, $G$-Brownian motion and stochastic differential equations driven by $G$-Brownian
motion ($G$-SDEs, in short).

 To establish Wang's Harnack inequality
using coupling by change of measures in the linear probability setting, the Girsanov transform plays a crucial role. In \cite{Hu, Osuka, Xu},   the Girsanov's theorem has been extended to the $G$-framework, and  the  Girsanov's formula has been derived
for  $G$-Brownian motion.
Recently,
  Hu et al. \cite{HWZ} studied the invariant and ergodic nonlinear expectations for $G$-diffusion processes.

In this paper, we investigate Wang's Harnack inequality and applications for the following
  $G$-SDE
\begin{equation}\label{SDE1}
dX_t=b(X_t) {d}\langle B \rangle_t+{ d} B_t,
\end{equation}
where $B_t$ is a $G$-Brownian motion, and $\langle B \rangle_t$ is the quadratic variation process
associated with $B_t$. Moreover,  we study  shift Harnack inequality and applications for the following
  $G$-SDE
\begin{equation}\label{SDE2}
dX_t=b(X_t) { d}t+{d} B_t,
\end{equation}where $B_t$ is a $G$-Brownian motion,

The paper is organized as follows. In Section 2, we recall some preliminaries on $G$-Brownian motion,  related stochastic calculus and transformation for $G$-Expection. In Section 3,
Wang's Harnack inequality and  shift Harnack inequality are established for the nonlinear Markov operator associated with \eqref{SDE1}  and  \eqref{SDE2} respectively.  As applications,
 the  $\sup$-kernel and invariant linear expectation for nonlinear Markov operator are investigated.

\section{Preliminaries}

\subsection{Sublinear expectation spaces}
In this section, we propose some preliminaries and notations which appeared in Peng \cite{peng2, peng4}.
Let $\Omega$ be a given set and $\cal{H}$ be a vector lattice of real valued functions defined
on $\Omega$, namely $c\in \cal{H}$ for each constant $c$ and $|X|\in \cal{H}$  if $X \in \cal{H}$. $\cal{H}$ is considered as the space of random variables.
\begin{definition} ({\rm  Sublinear expectation space})
 A sublinear expectation $\overline{E}$ on $\cal{H}$  is a functional $\overline{E}: \cal{H} \rightarrow \mathbb{R}$ satisfying the following properties: for all $X, Y \in  \cal{H},$  it holds that
 \begin{enumerate}
    \item [\rm (a)] Monotonicity: if $X \geq Y $ then $\overline{E}[X] \geq \overline{E}[Y ]$,
    \item [\rm (b)] Constant preservation: $\overline{E}[c] = c$,
    \item [\rm (c)] Sub-additivity: $\overline{E} [X + Y ] \leq \overline{E}[X] + \overline{E}[Y ]$,
    \item [\rm (d)] Positive homogeneity: $\overline{E} [\lambda X] = \lambda \overline{E} [X]$  for each $\lambda\geq 0$.
  \end{enumerate}
 Then,
 $(\Omega,\mathcal{H},\overline{E})$ is called  a sublinear expectation space.
\end{definition}
Let $\mathbb{S}^d$ be the collection of all $d\times d$ symmetric matrices,  $X$ be a $G$-normal distributed random vector, and
$G : \mathbb{S}^d \to R$ is defined by
$$G(A):=\frac{1}{2}\overline{E}[\langle AX,X\rangle]=\sup_{\gamma \in \Theta}\frac{1}{2}tr[\gamma \gamma^{\ast}A], \ A\in \mathbb{S}^d.
$$
Then the distribution of $X$ is characterized by
$$
u(t,x)=\overline{E}[\varphi(x+\sqrt{t}X)],\ \varphi\in C_{l,Lip}(\mathbb{R}^d),
$$
where  $C_{l,lip}(\mathbb{R}^n)$ be the space of all real functions $\varphi$ on $\mathbb{R}^{n}$ satisfying
\begin{eqnarray*}
|\varphi(x)-\varphi(y)|\leq C(1+|x|^{m}+|y|^{m})|x-y|, \ \ x,y\in \mathbb{R}^{n},
\end{eqnarray*}
for some $C>0, m \in \mathbb{N}$ depending on $\varphi$.

In particular, $\overline{E}[\varphi(X)] = u(1,0)$, where $u$ is the solution of the following parabolic
PDE defined on $[0,\infty)\times\mathbb{R}^d$:
\[
\left\{
\begin{array}{ll}
\frac{\partial u}{\partial t}-G(\frac{\partial^2 u}{\partial x^2})=0,\\
u(0,x)=\varphi(x).
\end{array}
\right.
\]
This parabolic PDE  is called a $G$-heat equation.
\subsection{ $G$-expectation and  $G$-Brownian motion}
Let $\Omega=C_0^d(\mathbb{R}^+)$  the space of all $\mathbb{R}^d$-valued continuous paths $(\omega_t)_{t\in \mathbb{R}^+}$, with $\omega_0 = 0$, equipped with the distance
\begin{equation*}
\rho(\omega^1,\omega^2):=\sum_{i=1}^\infty 2^{-i}\left[(\max_{t\in[0,i]}|\omega_t^1-\omega_t^2|)\wedge1\right].
\end{equation*}
Consider the
canonical process $B_t (\omega) = \omega_t, t \in [0,\infty)$, for $\omega \in \Omega.$ For each $T\in[0,\infty)$, let $\Omega_{T} =\{\omega_{\cdot \wedge T}: \omega \in \Omega\},$
and set
 $$L_{ip}(\Omega_T) =\{ \varphi(B_{t_1\wedge T}, \cdot \cdot \cdot, B_{t_n\wedge T}):n\in \mathbb{N}, t_1,\cdot \cdot \cdot, t_n\in [0,\infty),\varphi \in C_{l,lip}(\mathbb{R}^{d\times n})\}.$$
 Then $L_{ip}(\Omega_t)\subseteq L_{ip}(\Omega_T)$,  $t\leq T.$ Set
\begin{equation*}
  L_{ip}(\Omega)=\cup_{n=1}^{\infty}L_{ip}(\Omega_n).
\end{equation*}
Let $(\xi_i)_{i=1}^\infty$  a sequence of $d$-dimensional random vectors on a sublinear expectation space $(\tilde{\Omega},\mathcal{\tilde{H}},\tilde{E})$ such that $\xi_i$ is $G$-normal distributed and $\xi_{i+1}$ is independent from $(\xi_1 ,\cdot \cdot\cdot ,\xi_i )$ for each $i = 1,2,\cdot \cdot\cdot.$
We now introduce a sublinear expectation $\overline{E}$ defined on $L_{ip}(\Omega)$ via the following procedure: for each $X \in L_{ip}(\Omega)$ with
\begin{equation*}
  X=\varphi(B_{t_1}-B_{t_0},B_{t_2}-B_{t_1},\cdot \cdot \cdot,B_{t_n}-B_{t_{n-1}}), \ \varphi \in C_{l,lip}(\mathbb{R}^{d\times n}),
\end{equation*}
let
\begin{equation*}
 \overline{E}[\varphi(B_{t_1}-B_{t_0},B_{t_2}-B_{t_1},\cdot \cdot \cdot,B_{t_n}-B_{t_{n-1}})]=\tilde{E}[\varphi(\sqrt{t_1-t_0}\xi_1,\sqrt{t_2-t_1}\xi_1,\cdot \cdot \cdot,\sqrt{t_n-t_{n-1}}\xi_n)].
\end{equation*}
\begin{remark}\rm \label{remark2.2}
Let $\langle B\rangle_t:=(\langle B^{i},B^{j}\rangle_t)_{1\leq i,j\leq d}, 0\leq t\leq T$ be the quadratic variation of $B_t$.
In the 1-dimensional case, it holds that $\underline{\sigma}^2(t-s)\leq\langle B\rangle_t -\langle B\rangle_s\leq\overline{\sigma}^2(t-s), \ 0\leq s \leq t \leq T.$
\end{remark}

\begin{definition}\rm ($G$-expectation and  $G$-Brownian motion)
The  sublinear expectation $ \overline{E}: L_{ip}(\Omega) \to \mathbb{R}$ defined through the above procedure is called a $G$-expectation. The corresponding canonical process
$(B_t)_{t\geq 0}$ on the sublinear expectation space $(\Omega,L_{ip}(\Omega),\overline{E})$ is called a $G$-Brownian motion.
\end{definition}
\begin{remark}\rm
Let  $L_G^{p}(\Omega_T$) (respectively $L_G^p(\Omega)$) be  the completion of $L_{ip}(\Omega_T)$ (respectively $\Omega_t$) under the norm $(\overline{E}[|\cdot|^p])^{\frac{1}{p}}$. Then
$\overline{E}[\cdot]$ can be continuously extends to a sublinear expectation on $(\Omega,L_{G}^1(\Omega))$, which is  still denoted by $\overline{E}$.
\end{remark}
Let
\begin{equation*}\label{equa11} M_G^{p,0}([0,T]):=
\Bigg\{\eta_t:=\sum_{j=0}^{N-1} \xi_{j} I_{[t_j, t_{j+1})}; ~\xi_{j}\in L_{G}^p(\Omega_{t_{j}}),
 N\in\mathbb{N},\ 0=t_0<t_1<\cdots <t_N=T \Bigg\}.
\end{equation*}
For  $p\geq1$, let   $M_G^{p}([0,T])$ be the completion of $M_G^{p,0}([0,T])$ under the following norm
 $$\|\eta\|_{M_G^{p}([0,T])}=\left[\overline{E}\left(\int_{0}^{T}|\eta_{t}|^{p}dt\right)\right]
 ^{\frac{1}{p}}.$$

\subsection{Capacity and Quasi-Sure Analysis for $G$-Brownian Paths}
Denis et al. \cite{15} proved that there exists a weakly compact family $\{E_\theta:\theta \in \Theta\}$ of expectations introduced by  probability measures$\{ P_\theta: \theta \in \Theta\}$ defined on $(\Omega, \mathcal{B}(\Omega))$ such that
$$\overline{E}[X]=\sup_{\theta\in\Theta}E_{\theta}[X], \ X\in L_{ip}(\Omega).$$
Then the associated Choquet capacity is given by
$$c(A)=\sup_{\theta\in\Theta}P_\theta(A), \ A\in \mathcal{B}(\Omega).$$
\begin{definition} \rm \label{definition 2.1} (quasi-surely)
 A set $A\in\mathcal{B}(\Omega)$ is called polar if $c(A)=0$ and a property   holds  quasi-surely (q.s.)
if it holds outside a polar set.
\end{definition}
\begin{remark}\rm
Let $X$ and $Y$ be two random variables, we say that $X$ is a version of $Y$, if $X = Y $ q.s..
\end{remark}

\subsection{Stopping times}
In the sequel, we  introduce stopping times under $G$-expectation framework.
\begin{definition}\rm \label{definition 2.8}
Let $\mathcal{F}_{t}:=\mathcal{B}(\Omega)$, a stopping time $\tau$ relative to the filtration $(\mathcal{F}_{t})$ is a map on $\Omega$ with values in $[0,T]$, such that
 $$\{\tau\leq t\}\in\mathcal{F}_{t},\ t\leq T.$$
\end{definition}
\begin{lemma}\label{lemma1}{\rm{(\cite{Li})}}
For each stopping time $\tau$ and  $\eta\in M_{G}^{p}([0,T])$, we have $I_{[0,\tau]}(\cdot)\eta\in M_{G}^{p}([0,T])$,
and
\begin{equation*}
\int_{0}^{t\wedge\tau}\eta_{s}dB_{s}
=\int_{0}^{t\wedge\tau}\eta_{s}I_{[0,\tau]}(s)dB_{s}.
\end{equation*}
\end{lemma}

\subsection{Transformation for $G$-Expection}
To introduce the Girsanov type theorem under $G$-framework presented in \cite{Osuka, Xu},  using the  $G$-capacity  of \cite{15},
we need the $G$-Novikov's condition.
For $ h \in(M^2_G ([0, T]))^d$, let
\begin{eqnarray*}
&&M_t:=\exp \left\{  \int_0^th_s \cdot {d}B _s-\frac{1}{2}\int_0^th_s\cdot ({ d}\langle B \rangle_sh_s)\right \},\\
 &&\widehat{B} _t:=B_t-\int_0^{t}({d}\langle B \rangle _s h_s),  \ t\in[0,T].
\end{eqnarray*}
Set $\hat{L}_{ip}(\Omega_T) :=\{ \varphi(\hat{B}_{t_1\wedge T}, \cdot \cdot \cdot, \hat{B}_{t_n\wedge T}):n\in \mathbb{N}, t_1,\cdot \cdot \cdot, t_n\in [0,\infty),\varphi \in C_{l,lip}(\mathbb{R}^{d\times n}\}$. Let $\hat{L}_G^1(\Omega)$ be the completion of $\hat{L}_{ip}(\Omega_T)$ under the norm $\hat{E}[|\cdot|]$, and extend $\hat{E}$ to a unique sublinear
expectation on $ \hat{L}_G^1(\Omega).$
\begin{lemma} \label {lema2.10} {\rm (\cite{Osuka,Xu})}
 If $ h \in(M^2_G ([0, T]))^d$ satisfies $G$-Novikov's condition, for some $\epsilon_0 > 0$,
 \begin{equation}\label{Novikov}
 \overline{E}\left[ \exp \left\{ \left( \frac{1}{2}+\epsilon_0 \right )\int_0^Th_s\cdot({d}\langle B \rangle_sh_s)\right \}  \right]<\infty,
\end{equation}
then the process $M$ is a symmetric $G$-martingale.
\end{lemma}
We propose  Girsanov's formula for $G$-Brownian motion  as follows.
\begin{lemma}  \label{lemma3}{\rm(\cite{Osuka}) ($G$-Girsanov's formula)}
Assume that there exists $\sigma_0 > 0$ such that
 $$
 \gamma\gamma^\ast \geq \sigma_0I_d \  \ \  for \  \  \ all  \ \ \gamma \in \Theta,
 $$
 and that $M$ is a symmetric $G$-martingale on $(\Omega, L^1_G (\Omega), \overline{E}).$ Define a sublinear expectation $\widehat{E}$ by
 $$\widehat{E}[X]:=\overline{E}[XM_T]  \ \ \ for \ \ \ X \in  \hat{L}_{ip}^1(\Omega).$$
 Then $\widehat{B}_t$  is a $G$-Brownian motion on the sublinear expectation space $(\Omega,  \hat{L}_G^1(\Omega), \widehat{E})$.
\end{lemma}
Moreover, Hu et al.  \cite{Hu} construct an auxiliary extended $\widetilde{G}$-expectation space $(\widetilde{\Omega}_T, L_{\widetilde{G}}^1, \overline{E}^{\widetilde{G}})$ with $\widetilde{\Omega}_T=C_0([0,T],\mathbb{R}^2)$ and
$$
\widetilde{G}(A)=\frac{1}{2}\sup_{\underline{\sigma}^2\leq v\leq\overline{\sigma}^2}tr\left[A
\left(
  \begin{array}{cc}
    v & 1 \\
    1 & v^{-1} \\
  \end{array}
\right)\right], \ A\in \mathbb{S}^{2}.
$$
Let $(B_t,\bar{B}_t)$ be the canonical process in the extended space. Then $\langle B_t, \bar{B}_t\rangle=t$,
and $$\overline{E}^{\widetilde{G}}[\xi]=\overline{E}[\xi],  \ \xi \in L_G^1(\Omega_T). $$
\begin{lemma} \label{lemma4}  {\rm(\cite{Hu})}
Let $(h_t)_{t\geq0}$  be a bounded process,
 then the process  $\widetilde{B}_t:=B_t+\int_0^th_sds$  is a $G$-Brownian motion under $\widetilde{E}$ for
$$\widetilde{E}[X]=\overline{E}^{\widetilde{G}}\left[X\exp{\left(  \int_0^Th_s{ d}\bar{B} _s-\frac{1}{2}\int_0^Th^2_s{d}\langle \bar{B} \rangle_s\right)}\right], \ X\in L_{G}^1(\Omega_T).$$
\end{lemma}
\begin{remark} \rm \label{remark2.8} We  should remark that the
$\bar{B}_t$ is a  $\hat{G}$-Brownian motion under $\overline{E}^{\widetilde{G}}$ with
$\hat{G}(A)=\frac{1}{2}\sup_{\overline{\sigma}^{-2}\leq v\leq\underline{\sigma}^{-2}}tr[Av]$, $A\in \mathbb{S}^{1}$.
\end{remark}

\section{Main Results}
For a family of probability measures $\{\mu_{x,\theta}:x\in \mathbb{R}^d, \theta \in \Theta\}$  on $(\mathbb{R},\mathcal{B}(\mathbb{R}^d))$,  define
\begin{equation*}
\bar{P}f(x)=\sup_{\theta \in \Theta}{P}_\theta f(x)=\sup_{\theta \in \Theta}\int_{\mathbb{R}^d} f(y)\mu_{x,\theta}(dy), \ f\in \mathcal{B}_b(\mathbb{R}^d),
\end{equation*}
where ${P}_\theta $ is a linear Markov operator.

 We aim to establish the following  Harnack-type inequality introduced by Feng-Yu Wang:
\begin{equation*}
\Phi(\bar{P}f(x))) \leq\bar{P}\Phi(f(y)) e^{\Psi(x,y)}, \  x,y \in \mathbb{R}^d, f\in \mathcal{B}_b^+(\mathbb{R}^d),
\end{equation*}
where $\Phi$ is a  nonnegative convex function on $[0,\infty)$ and $\Psi$ is a  nonnegative  function on $\mathbb{R}^d\times \mathbb{R}^d$.

In the   setting of  $G$-SDEs, we establish this type inequality for the associated Markov operator $\bar{P}_T$. For simplicity, we consider the 1-dimensional $G$-Brownian motion case,
but our results and methods still hold for the case $d > 1$. More precisely,
consider the following $G$-SDE
\begin{equation}\label{0}
dX_t=b(X_t) {d}\langle B \rangle_t+{ d} B_t,
\end{equation}
where $B_t$ is a $G$-Brownian motion,   $\langle B \rangle_t$ is the quadratic variation process
associated with $B_t$ for  $\overline{\sigma}^2 = \overline{E} [B^2_1 ] \geq -\overline{E} [-B^2_1 ] = \underline{\sigma}^2 > 0$, and $b:\mathbb{R} \to \mathbb{R}$ satisfies
\begin{itemize}
  \item [(H1)] $|b(x)-b(y)|\leq K|x-y|$, $  \ x, y \in \mathbb{R}$
\end{itemize}
for some constant $K>0$.
From \cite{peng1}, under (H1) the $G$-SDE \eqref{0} has a unique solution for any initial value.

In what follows,  for $T>0$, we define
\begin{equation}\label{deP}
\bar{P}_Tf(x)=\overline{E}f(X_T^x), \  f\in \mathcal{B}_b^+(\mathbb{R}),
\end{equation}
where  $ X_T^x$ solves \eqref{0} with initial value $x$.
We have  the following result.
\begin{theorem} \label {th1}
Under {\rm(H1)},  for any nonnegative $ f \in \mathcal{B}_b^+(\mathbb{R})$ and $ T> 0, x, y \in  \mathbb{R},$  it holds that
\begin{equation} \label {th1eq1}
(\bar{P}_Tf)^p(y)\leq \bar{P}_Tf^p(x)\exp \left\{\frac{pK\overline{\sigma}^4(1-e^{-2\underline{\sigma}^2KT})}{(p-1)^2\underline{\sigma}^6(1-e^{-2\overline{\sigma}^2KT})^2}|x-y|^2 \right\}.
\end{equation}
\end{theorem}
\textbf{Proof.} We use the coupling by change of measures as explained in \cite{W1},
consider the following coupled stochastic differential equations
\begin{eqnarray}\label{1}
 dY_t&=&b(Y_t) {d}\langle B \rangle _t+{d} B_t+u_t{d} \langle B \rangle _t, \ Y_0=y, \\
dX_t&=&b(X_t) {d}\langle B \rangle _t+{d} B_t, \ X_0=x, \nonumber
\end{eqnarray}
where $u_t=\eta_t\cdot\frac{X_t-Y_t}{|X_t-Y_t|}1_{t\leq\tau}$, $\tau$ is the coupling time of $X$ and $Y$ defined by
$$\tau=\inf\{t\geq 0:X_t=Y_t  \},$$
and $u_\cdot$ is a force can make the two processes $X$ and $Y$ move together before time $T$ .

From assumption (H1)  and the expression of  $\eta_t$,  \eqref{1}  has a unique solution.

By (H1), we have
 $${ d}|X_t-Y_t|\leq K|X_t- Y_t| { d}\langle B \rangle _t- \eta_t{d}\langle B \rangle _t, \ t< \tau.$$
Then
\begin{equation*}
e^{-K\langle B\rangle _{T\wedge \tau}}|X_{T\wedge \tau}- Y_{T\wedge \tau}|\leq |x-y|-\int_0^{T\wedge \tau}e^{-K\langle B\rangle _{t}}\eta_t{ d}\langle B\rangle _t.
\end{equation*}
Taking
\begin{equation*}
\eta_t=\frac{|x-y|e^{-K\langle B\rangle _t}}{\int_0^Te^{-2K\langle B\rangle _t}{ d}\langle B \rangle _t}, \   \   \    t\in[0,T],
\end{equation*}
we have
\begin{equation*}
|x-y|-\int_0^{T}e^{-K\langle B \rangle _t}\eta_t{d}\langle B \rangle _t=0,
\end{equation*}
which implies $\tau\leq T$, and thus $X_T=Y_T$.

By Remark \eqref{remark2.2}, we have
\begin{eqnarray}\label{Har3}
&&\exp \left\{\int_0^{T} |u_s|^2{d}\langle B \rangle_s\right\} \nonumber \\
&=&\exp \left\{ \overline{\sigma}^2\int_0^{T} \left|\frac{|x-y|e^{-K\langle B \rangle_s}}{\int_0^Te^{-2K\langle B \rangle_s}{d}\langle B \rangle_s}\right|^2{ d}s\right\} \nonumber \\
&\leq&\exp \left\{ \frac{\overline{\sigma}^2}{\underline{\sigma}^4} \frac{\int_0^{T}\left||x-y|e^{-K\langle B \rangle_s}\right|^2{ d}s}{\left|\int_0^Te^{-2K\langle B \rangle_s}{ d}s\right|^2}\right\} \nonumber \\
&\leq&\exp \left\{\frac{2K\overline{\sigma}^4(1-e^{-2\underline{\sigma}^2KT})}{\underline{\sigma}^6(1-e^{-2\overline{\sigma}^2KT})^2}|x-y|^2 \right\}.
\end{eqnarray}
Letting $\epsilon_0=\frac{1}{2}$, we have
\begin{equation*}
 \overline{E}\left[ \exp \left\{ \int_0^{T}  \left|u_s \right|^2{d} \langle B \rangle_s\right\}\right]
\leq \exp \left\{\frac{2K\overline{\sigma}^4(1-e^{-2\underline{\sigma}^2KT})}{\underline{\sigma}^6(1-e^{-2\overline{\sigma}^2KT})^2}|x-y|^2\right\}<\infty,
\end{equation*}
which satisfies $G$-Novikov's condition.

 Let
 \begin{equation*}
M_T=\exp \left\{-\int_0^Tu_s{d}B_s -\frac{1}{2}\int_0^T |u_s|^2{d}\langle B \rangle_s\right\}.
\end{equation*}
Define a sublinear expectation $\widehat{E}$ by
$\widehat{E}[\xi]:=\overline{E}[\xi M_T]$.
 By Lemma \ref{lema2.10} and Lemma \ref{lemma3}, the process
\begin{equation*}
 \widehat{B} _t:=B_t+\int_0^{t}u_s{d}\langle B \rangle _s, \ \ \ t\geq 0
\end{equation*}
is a $G$-Brownian motion under $\widehat{E}$.

Moreover, Girsanov's formula also implies that
$$\langle \widehat{B}\rangle_{\widehat{E}}=\langle B \rangle_{\overline{E}}. $$
Then, $Y_t$ can be reformulated by
\begin{equation*}
{ d}Y_t=b(Y_t) { d}\langle \widehat{B}\rangle_t+{ d} \widehat{B}_t.
\end{equation*}
So, $\bar{P}_Tf(y)=\overline{E}f(X_T^y)=\widehat{E}f(Y_T^y)=\widehat{E}f(X_T^x)=\overline{E}(M_Tf(X_T^x))$.

Using H\"{o}lder's inequality, we have
\begin{eqnarray}\label{Har1}
(\bar{P}_Tf)^p(y)&=& (\overline{E}[M_Tf(X_T^x)])^p\nonumber\\
&\leq& (\overline{E}[f^p(X_T^x)])\left(\overline{E}\left[M_T^{\frac{p}{p-1}}\right]\right)^{\frac{p-1}{p}}.
\end{eqnarray}
Now we  estimate the moment of $M_T$. It holds that
\begin{eqnarray}\label{Har2}
\overline{E}\left[M_T^{\frac{p}{p-1}}\right]&=&\overline{E}\exp \left\{-{\frac{p}{p-1}}\int_0^{T} u_s{d}B_s -\frac{p}{2(p-1)}\int_0^{T} |u_s|^2{ d}\langle B \rangle_s\right\} \nonumber \\
&=&\overline{E}\exp \Bigg\{-{\frac{p}{p-1}}\int_0^{T} u_s{d}B_s -\frac{p^2}{2(p-1)^2}\int_0^{T} |u_s|^2{ d}\langle B \rangle_s\nonumber \\
&&+\frac{p}{2(p-1)^2}\int_0^{T} |u_s|^2{ d}\langle B \rangle_s\Bigg\}.
 \end{eqnarray}
From \eqref{Har3}, we have
\begin{equation*}
\exp \left\{ \frac{p}{2(p-1)^2}\int_0^{T} |\eta_s|^2{d}\langle B \rangle_s\right\}
\leq\exp \left\{\frac{pK\overline{\sigma}^4(1-e^{-2\underline{\sigma}^2KT})}{(p-1)^2\underline{\sigma}^6(1-e^{-2\overline{\sigma}^2KT})^2}|x-y|^2 \right\}.
\end{equation*}
Substituting this into \eqref{Har2}, we have
\begin{equation}\label{Har4}
\overline{E}\left[M_T^{\frac{p}{p-1}}\right]\leq\exp \left\{\frac{pK\overline{\sigma}^4(1-e^{-2\underline{\sigma}^2KT})}{(p-1)^2\underline{\sigma}^6(1-e^{-2\overline{\sigma}^2KT})^2}|x-y|^2 \right\}.
 \end{equation}
Combing \eqref{Har1} and  \eqref{Har4}, we prove \eqref{th1eq1}.

To obtain the shift Harnack inequality, we consider the  following $G$-SDE
\begin{equation}\label{SDE1'}
dX_t=b(X_t) {\rm d}t+{d} B_t,\ X_0=x,
\end{equation}
where  $B_t$ is a $G$-Brownian motion, and $b:\mathbb{R}\to \mathbb{R}$ satisfies
the assumption (H1).
Then, the $G$-SDE \eqref{SDE1'} has a unique solution for any initial value.

For $T>0$, we define
\begin{equation*}
\bar{P}_Tf(x)=\overline{E}f(X_T^x), \  f\in \mathcal{B}_b^+(\mathbb{R}),
\end{equation*}
where  $ X_T^x$ solves \eqref{SDE1'} with initial value $x$. By the definition of $\overline{E}^{\widetilde{G}}$ (in section 2.6), we have
$$\overline{E}[f(X_T^x)]=\overline{E}^{\widetilde{G}}[f(X_T^x)]=:\bar{P}_T^{\widetilde{G}}f(x).$$
We have the following result.
\begin{theorem}
Under {\rm(H1)},  for any nonnegative $ f \in  \mathcal{B}_b^+(\mathbb{R})$ and $ T > 0, x, y \in   \mathbb{R},$  the following shift Harnack inequality holds
\begin{equation}\label{th2eq}
(\bar{P}_Tf(x))^p \leq(\bar{P}_Tf^p(v+\cdot))(x)\exp \left\{ \frac{pv^2}{2\underline{\sigma}^2(p-1)}\left(\frac{1}{T}+K+\frac{K^2T}{3}\right)\right\}.
\end{equation}
\end{theorem}
\textbf{Proof.}
Let $Y_t=X_t+\frac{t}{T}v$ with  $Y_0=X_0=x$ and
\begin{equation*}
R_t=\exp{\left(-\int_0^t\frac{v}{T}+b(X_s)- b(Y_s){ d}\bar{B} _s-\frac{1}{2}\int_0^t(\frac{v}{T}+b(X_s)- b(Y_s))^2{d}\langle \bar{B} \rangle_s\right)},
\end{equation*}
where $\bar{B} _t$ is $G$-Brownian motion under $\overline{E}^{\widetilde{G}}$, which is an auxiliary process, one can see in section 2.6 for details.

From (H2), we have
\begin{equation}\label{bounded}
\left|\frac{v}{T}+b(X_s)- b(Y_s)\right|\leq \frac{1+Ks}{T}|v|.
\end{equation}
By Lemma\ref{lemma4},
\begin{equation*}
 \widetilde{B}_t:=B _t+\int_0^t  \left(\frac{v}{T}+b(X_s)- b(Y_s)\right)  {d} s
\end{equation*}
is a $G$-Brownian motion under $\widetilde{E}$ with
$$\widetilde{E}[\xi]=\overline{E}^{\widetilde{G}}\left[\xi R_T\right],  \ \xi \in L_{G}^1(\Omega_T).$$
Then
\begin{equation*}
 {\rm d}Y_t=b(Y_t) { d}t+{ d}\widetilde{ B} _t,   \   Y_0=x.
\end{equation*}
That is  $Y_T=X_T+v$ under  $\widetilde{E}$.

Then for  $ f \in \mathcal{B}_b^+(\mathbb{R})$, $p\geq1$,  by H\"{o}lder inequality, we have
\begin{eqnarray}\label{h1}
\left(\bar{P}_Tf\right)^p(x)
&=& \left(\widetilde{E}[f(Y_T^x)]\right)^p\nonumber\\
&=& \left(\widetilde{E}[f(X_T^x+v)]\right)^p\nonumber\\
&=& \left(\overline{E}^{\widetilde{G}}[R_Tf(X_T^x+v)]\right)^p\nonumber\\
&\leq&\left(\bar{P}_T^{\widetilde{G}}f^p(v+\cdot)\right)(x)\left(\overline{E}^{\widetilde{G}}[R_T^{\frac{p}{p-1}}]\right)^{p-1}\nonumber\\
&=&\left(\bar{P}_Tf^p(v+\cdot)\right)(x)
\left(\overline{E}^{\widetilde{G}}\left[R_T^{\frac{p}{p-1}}\right]\right)^{p-1}.
\end{eqnarray}
Letting $h_s=\frac{v}{T}+b(X_s)- b(Y_s)$, by Remark \ref{remark2.8} and \eqref{bounded},  we have
\begin{eqnarray*}
\exp \left\{ \frac{p}{2(p-1)^2}\int_0^{T} |h_s|^2{ d}\langle \bar{B} \rangle_s\right\}
&\leq&\exp \left\{ \frac{p\underline{\sigma}^{-2}}{2(p-1)^2}\int_0^{T} |h_s|^2{ d}s\right\} \nonumber \\
&=&\exp \left\{ \frac{p\underline{\sigma}^{-2}v^2}{2(p-1)^2}\left(\frac{1}{T}+K+\frac{K^2T}{3}\right)\right\}.
\end{eqnarray*}
Similarly with the proof of Theorem \ref{th1},  we have
 \begin{equation}\label{h3}
\overline{E}^{\widetilde{G}}\left[R_T^{\frac{p}{p-1}}\right]
\leq\exp \left\{ \frac{p\underline{\sigma}^{-2}v^2}{2(p-1)^2}\left(\frac{1}{T}+K+\frac{K^2T}{3}\right)\right\}.
\end{equation}
It follows from \eqref{h1} and \eqref{h3} that \eqref{th2eq} holds.

\subsection{Applications of Harnack  and shift  Harnack Inequalities }
In this subsection, we give some applications of Harnack  and shift  Harnack inequalities for  invariant linear  expectation and  $\sup$-kernel estimates. Before that,
due to technical difficulties, we need the following invariant linear  expectation, let us define
the  (quasi) invariant linear (nonlinear) expectation and $\sup$-kernel  of the  operator $\bar{P}$.
\begin{definition}\rm \label{def4.1}
Let $E$ be a linear (nonlinear) expectation, 
and  $\bar{P}$ be a nonlinear operator defined on $\mathcal{B}_b^+(\mathbb{R}^d)$.
 \begin{enumerate}
    \item [\rm (1)]  $E$ is called a quasi-invariant linear (nonlinear) expectation of $\bar{P}$, if there exists a function $0\leq g \in \mathcal{B}_b(\mathbb{R}^d)$ with $E[g]<\infty$, such that
     \begin{equation*}
    E[\bar{P}f]\leq E[gf], \  \ 0\leq f \in \mathcal{B}_b(\mathbb{R}^d).
     \end{equation*}
    Moreover, if
      \begin{equation*}
E[(\bar{P}f)]= {E}[f],  \  \ 0\leq f \in \mathcal{B}_b(\mathbb{R}^d),
     \end{equation*}
     then $E$ is called an  invariant   linear (nonlinear) expectation of $\bar{P}$.
     \item [\rm (2)] A function $p$ on $\mathbb{R}^d\times \mathbb{R}^d$ is    called   the $\sup$-kernel (or $\sup$-density) of $\bar{P}$ with respect to $ {E}$, if
     \begin{equation*}
    \bar{P}f(x)\leq E[p(x,\cdot) f(\cdot)], \  \ 0\leq f \in \mathcal{B}_b(\mathbb{R}^d), \ x\in \mathbb{R}^d.
     \end{equation*}
\end{enumerate}
\end{definition}
To illustrate the above definition, we consider an example as follows.
\begin{example}\rm We consider the following Ornstein-Uhlenbeck process driven by $G$-Brownian motion:
for each $x\in \mathbb{R}^d$,
\begin{equation*}
Y_t^x=x-\alpha\int_0^t Y_s^xds+B_t, \ t\geq 0,
\end{equation*}
where $B_t$ is a $d$-dimensional  $G$-Brownian motion.
 Hu et al. \cite{HWZ} proved the unique invariant expectation for $G$-Ornstein-Uhlenbeck process $Y$ is the
$G$-normal distribution of $\sqrt{\frac{1}{2\alpha}}B_1$.
\end{example}

\begin{example}\rm \label{ex3.4}
For $\theta \in \left[\frac{1}{2},1\right],$ let $W_t$ be the stand 1-dimensional Brownian motion, consider
 the following SDE,
\begin{equation*}
dX_t=-\theta X_tdt+\sqrt{2\theta}dW_t,\ X_0=x.
\end{equation*}
Then, $X_t=e^{-\theta t}x+\int_0^t\sqrt{2\theta}e^{-\theta (t-s)}dW_s$, $X_t \to N(0,1)$ in distribution as $t\to \infty.$

Let
\begin{equation*}
P_\theta f(x)=\int_R\frac{1}{\sqrt{2\pi(1-e^{-2\theta })}}\exp\left\{-{\frac{{(z-e^{\theta }x)^2}}{2(1-e^{-2\theta })}}\right\}f(z)dz, \ \theta \in \left[\frac{1}{2},1\right].
\end{equation*}
Therefore,
\begin{eqnarray*}
P_\theta f(x)&=&E_0\left[\frac{\frac{1}{\sqrt{2\pi(1-e^{-2\theta })}}\exp\left\{-{\frac{{(\cdot-e^{\theta }x)^2}}{2(1-e^{-2\theta })}}\right\}f(\cdot)}
{\frac{1}{\sqrt{2\pi}}e^{-\frac{(\cdot)^2}{2}}}\right]\\
&=&E_0\left[\frac{\frac{1}{\sqrt{1-e^{-2\theta }}}\exp\left\{-{\frac{{(\cdot-e^{\theta }x)^2}}{2(1-e^{-2\theta })}}\right\}f(\cdot)}
{e^{-\frac{(\cdot)^2}{2}}}\right],
\end{eqnarray*}
where $E_0[f]=\int_{\mathbb{R}}fd\mu_0, \mu_0=N(0,1).$

Moreover,
\begin{equation*}
\frac{\frac{1}{\sqrt{1-e^{-2\theta }}}\exp\left\{-{\frac{{(z-e^{\theta }x)^2}}{2(1-e^{-2\theta })}}\right\}}
{e^{-\frac{z^2}{2}}}
\leq\frac{e^{\frac{z^2}{2}}}{\sqrt{(1-e^{- 1})}}.
\end{equation*}
Then, $p(x,y)=\frac{e^{\frac{y^2}{2}}}{\sqrt{1-e^{- 1}}}$ is a $\sup$-kernel of $\bar{P}$ with respect to $ {E_0}$, where $\bar{P}f=\sup_{\theta \in \Theta} P_\theta f$.
\end{example}

\subsubsection{Applications of Harnack Inequalities}
Now, we consider following  Harnack-type inequality
\begin{equation} \label {ap1}
\Phi(\bar{P}f(x)) \leq\bar{P}\Phi(f(y)) e^{\Psi(x,y)},
\end{equation}
where $\Phi$ is a  nonnegative convex function on $[0,\infty)$ and $\Psi$ is a  nonnegative  function on $\mathbb{R}^d\times \mathbb{R}^d$.

\begin{theorem}\label{proposition3} Let $E$ be  a  quasi-invariant linear expectation of $\bar{P}$, and  $\Phi \in C^1([0,\infty))$ be an increasing function with $\Phi'(1) > 0$ and $\Phi(\infty) := \lim_{r \to \infty }\Phi(r) = \infty $ such that \eqref{ap1} holds.
 \begin{enumerate}
     \item [\rm (1)]  If  $\lim_{y\rightarrow x}\{\psi(x,y)+\psi(y,x)\}=0 $ holds for all $x \in \mathbb{R}^d$, then $\bar{P}$ is strong Feller, i.e. $\bar{P}\mathcal{B}_{b}(\mathbb{R}^d)\subset C_b(\mathbb{R}^d)$.
  \item [\rm (2)] Let $\bar{P}f(x)=\sup_{\theta \in \Theta}P_\theta f(x)$. Then for all  $ \theta \in \Theta,$   $P_\theta$ has  a kernel $p_\theta$ with respect to ${E}$,
    $\bar{P}$ has a $\sup$-kernel $p$ with respect to $ {E}$, and  every invariant linear  expectation
of $\bar{P}$ is absolutely continuous with respect to ${E}$.
 \item [\rm (3)]If there exists $K>0$ such  that $\frac{1}{K}p_{\theta_1}(x,y)\leqslant p_{\theta_2}(x,y)\leqslant K p_{\theta_1}(x,y)$,   $ \theta_1,  \theta_2 \in \Theta$, $x,y \in \mathbb{R}^d$, where $p_{\theta}(x,y)$ is defined  in (2), then
 $\bar{P}$ has at most one  invariant linear  expectation, and if it has one, a $\sup$-kernel of
$\bar{P}$ with respect to the invariant linear expectation is strictly positive.
 \item [\rm (4)] If $r\Phi^{-1}(r)$ is convex for $r\geqslant0$, then a $\sup$-kerner $p$ of $\bar{P} $  with respect to $E$ satisfies
 \begin{equation*}
E\left[ p(x,\cdot) p(y,\cdot)\right]\geq e^{-\Psi(x,y)}.
\end{equation*}
  \item [\rm (5)] If ${E}$ is a    invariant  linear expectation of $\bar{P}$, then
   \begin{equation*}
\sup_{f \in \mathcal{B}^+_{b}(\mathbb{R}^d), E[\Phi(f)]\leq1} \Phi(\bar{P}f(x))\leq \frac{1}{E[e^{-\Psi(x,\cdot)}]}.
\end{equation*}
\end{enumerate}
\end{theorem}
\textbf{Proof.}
\begin{enumerate}

\item [(1)]  Let $0<f \in \mathcal{B}_{b}(\mathbb{R}^d)$. Applying \eqref{ap1} to $f=1+\epsilon f$ for $\epsilon >0$, we have
\begin{equation*}
\Phi(1+\epsilon\bar{P}f(x)) \leq \{\bar{P}(\Phi(1+\epsilon f(y)))\}e^{\Psi(x,y)}, \ x,y \in \mathbb{R}^d, \ \epsilon >0.
\end{equation*}
By a Taylor expansion, we get
\begin{equation}\label{lian}
\Phi(1)+\epsilon\Phi'(1)\bar{P}f(x)+o(\epsilon) \leq \{\bar{P}(\Phi(1)+\epsilon \Phi'(1)f(y)+o(\epsilon))\}e^{\Psi(x,y)}
\end{equation}
for small $\epsilon > 0$. Letting $y\to x$, we have
\begin{equation*}
\epsilon\bar{P}f(x) \leq \epsilon \liminf_{y\rightarrow x}\bar{P}f(y)+o(\epsilon).
\end{equation*}
Then $\bar{P}f(x) \leq  \liminf_{y\rightarrow x}\bar{P}f(y)$  for all   $x\in \mathbb{R}^d$.

Moreover,  letting $x\to y$ in \eqref{lian}, we have
 $\bar{P}f(y) \geq \limsup_{x\rightarrow y}\bar{P}f(x)$  for all   $ y\in \mathbb{R}^d$.
Consequently, $\bar{P} f$ is continuous.

\item [(2)]  To prove the existence of a $\sup$-kernel, 
 it suffices to prove  $\bar{P}1_A(x)\leq  {E}[p(x,\cdot) 1_A]$ for some positive function $p$ on $\mathbb{R}^d\times \mathbb{R}^d$.

 We firstly prove that ${E}[1_A]=0$  implies $\bar{P}1_A\equiv0$. Applying \eqref{ap1} to $f=1+n1_A$, we have
 \begin{equation*}
\Phi(1+n\bar{P}1_A(x)) \leq \bar{P}\Phi(1+n1_A(y)) e^{\Psi(x,y)}, \ x,y \in \mathbb{R}^d.
\end{equation*}
Then
 \begin{equation*}
\Phi(1+n\bar{P}1_A(x))e^{-\Psi(x,y)} \leq \bar{P}\Phi(1+n1_A(y)), \ x,y \in \mathbb{R}^d.
\end{equation*}
 It follows from   ${E}$ is   a quasi-invariant  expectation of $\bar{P}$ and $E[g]<\infty$ that
 \begin{eqnarray*}
\Phi(1+n\bar{P}1_A(x)) &\leq& \frac{{E}[\bar{P}\Phi(1+n1_A(\cdot))]}{{E}[e^{-\Psi(x,\cdot)}] }\\
&=&\frac{{E}[\Phi(1+n1_A(\cdot))g]}{{E}[e^{-\Psi(x,\cdot)}] }\\
&\leq&\frac{{E}[g\Phi(1)]}{{E}[e^{-\Psi(x,\cdot)}] }\\
&=&\frac{\Phi(1){E}[g]}{{E}[e^{-\Psi(x,\cdot)}] }\\
&{<}&\infty.
\end{eqnarray*}
 Since $ \lim_{n \to \infty }\Phi(1+n) = \infty $, which  implies that $\bar{P}1_A(x) = 0$, $x \in \mathbb{R}^d$.

As $\bar{P}1_A(x) = \sup_{\theta \in \Theta}P_\theta 1_A(x)$, then $P_\theta 1_A(x)=0$,
 for all   $\theta  \in \Theta$. Therefore,
 there exists a $p_\theta$ on $\mathbb{R}^d\times \mathbb{R}^d$, such that ${P}_\theta f(x)= {E}[p_\theta(x,\cdot) f(\cdot)]$, $x\in \mathbb{R}^d$.

 Moreover,
 \begin{equation*}
\bar{P}f(x)
 =\sup_{\theta \in \Theta}P_\theta f(x)
 =\sup_{\theta \in \Theta} {E}[p_\theta(x,\cdot) f(\cdot)]
 \leq  {E}[p(x,\cdot) f(\cdot)], \   0<f \in \mathcal{B}_{b}(\mathbb{R}^d),
 \end{equation*}
 where $p(x,\cdot)=\sup_{\theta \in \Theta}p_\theta(x,\cdot)$, we assume  ${E}[p(x,\cdot) ]<\infty$.

 Furthermore, for any invariant  expectation ${E}_0$ of $\bar{P}$, if ${E}[1_A] = 0$, then $\bar{P}1_A= 0$ implies that
  ${E}_0[1_A]  = {E}_0[\bar{P}1_A] = 0$. Therefore, ${E}_0$  is absolutely continuous with respect to ${E}$.

 \item [(3)]
   Let $p$ be a  $\sup$-kernel of $\bar{P}$ with respect to every invariant linear expectation  ${E}_0$.
  \begin{equation*}
\bar{P}f(x) \leq  {E}_0[p(x,\cdot) f(\cdot)]=:\widetilde{P}f(x), \  0<f \in \mathcal{B}_{b}(\mathbb{R}^d).
 \end{equation*}
We aim wo prove $p>0$.
   In fact, from the definition of
    $\widetilde{P}$, then ${E}_0[1_A]  = 0$ implies that
$\widetilde{P}1_A= 0$.  To this end,
    it suffices to show that for any $x\in \mathbb{R}$,   $\widetilde{P}1_{A}(x) = 0$ implies that  ${E}_0[1_{A}(x)] = 0$.
   Since $\widetilde{P}1_{A}(x)  \geq \bar{P}1_A(x)$, it suffices to show that   $\bar{P}1_{A}(x) = 0$ implies that  ${E}_0[1_{A}(x)] = 0$.
    Since $\bar{P}1_A(x) = 0,$ by applying \eqref {ap1}  to $f = 1+n1_{A}$, we obtain
    \begin{equation*}
\Phi(1+n\bar{P}1_{A}(y))\leq \bar{P}\Phi(1+n1_{A})(x) e^{\Psi(y,x)}=\Phi(1) e^{\Psi(y,x)}.
\end{equation*}
  Letting $n \to \infty $, we conclude that $\bar{P}1_A\equiv 0$,   then   ${E}_0[1_{A}]= {E}_0[\bar{P}1_{A}] = 0$, which implies the $\sup$-kernel $p(x,y)>0.$

Next we prove the uniqueness of invariant expectation.
    Let  ${E}_1$ is a  another  invariant  linear expectation. From (2),  there exists a function $f$, for any $g$, such that
     ${E}_1[g]={E}_0[fg]$, and $ E_0[f]=1$. We aim to prove that $f=1,  {E}_0-a.e.$.
    Let  $p(x,y)>0$ be a $\sup$-kernel of $\bar{P}$ with respect to  ${E}_0$. Then, $\bar{P}f(x)\leq {E}_0[p(x,\cdot)f(\cdot)]$.

      Moreover, since $\bar{P}f(x)= \sup_{\theta \in \Theta}P_\theta f(x)$, by (2), we know that ${E_0}[1_A]=0$  implies ${P_\theta}1_A=0$.

       Therefore, for any  $\theta \in \Theta$,
   we have
   \begin{equation*}
   \bar{P} f(x)\geq {P_\theta} f(x)={E}_0[p_\theta(x,\cdot)f(\cdot)].
   \end{equation*}
    Let $P^{\ast}_{\theta}f(x)=\int_{\mathbb{R}^d}f(y)P_\theta^{\ast}(x,dy)={E}_0[p_\theta(\cdot,x)f(\cdot)]$, $x\in \mathbb{R}^d$.

    For any $0<h\in L^1(E_0), 0<g \in \mathcal{B}_{b}(\mathbb{R}^d)$, by Fubini theorem, we have
    \begin{equation}\label{Fubini}
    {E}_0 [P^{\ast}_{\theta}hg]={E}_0[{E}_0[p_\theta(\cdot,x)h(\cdot)]g(x)]={E}_0[ {E}_0 [p_\theta(y,\cdot)g(\cdot)]h(y)]={E}_0[ hP_\theta g].
    \end{equation}
 Since ${E}_0$ is $\bar{P}$-invariant, and $P_{\theta}1=1,$  by \eqref{Fubini}, we have
    \begin{equation*}
{E}_0[(P^{\ast}_{\theta}1)g]={E}_0[P_{\theta}g]\leq{E}_0[\bar{P}g]={E}_0[g], \  0<g \in \mathcal{B}_{b}(\mathbb{R}^d),
  \end{equation*}
which implies $P^{\ast}_{\theta}1\leq 1$, ${E}_0$-a.e..

As $f\in L^1(E_0),$ by \eqref{Fubini}, we have
 \begin{equation*} \label{E3}
{E}_0[(P^{\ast}_{\theta}f)g]={E}_0[fP_{\theta}g]={E}_1[P_{\theta}g]\leq {E}_1[\bar{P}g]={E}_1[g]={E}_0[fg],  \ 0<g \in \mathcal{B}_{b}(\mathbb{R}^d).
\end{equation*}

which means $P^{\ast}_{\theta}f\leq f$, ${E}_0$-a.e..

By $P^{\ast}_{\theta}1\leq 1$, ${E}_0$-a.e., and H\"older inequality, we have
 \begin{equation} \label{E0}
 P^{\ast}_{\theta} \sqrt{f+1}\leq \sqrt{P^{\ast}_{\theta}(f+1)} \sqrt{P^{\ast}_{\theta}1}
 \leq \sqrt{P^{\ast}_{\theta}f+P^{\ast}_{\theta}1}\leq \sqrt{P^{\ast}_{\theta}f+1}, \ {E}_0\rm{-a.e.}.
\end{equation}

 Furthermore,  by \eqref{Fubini} and  $P^{\ast}_{\theta}f\leq f$, ${E}_0$-a.e.,  we obtain
\begin{equation} \label{E1}
{E}_0\left[P^{\ast}_{\theta}\sqrt{f+1}\right]={E}_0\left[\sqrt{f+1}P_{\theta}1\right]={E}_0\left[\sqrt{f+1}\right]\geq{E}_0\left[\sqrt{P^\ast_{\theta}f+1} \right].
\end{equation}
From \eqref{E0} and  \eqref{E1}, we get
\begin{equation} \label{E2}
P^{\ast}_{\theta}\sqrt{f+1}= \sqrt{P^\ast_{\theta}f+1}, \ {E}_0 \rm{-a.e.}.
\end{equation}
By recalling  \eqref{E0} again, we have
  \begin{equation}\label{holdereq}
 P^{\ast}_{\theta} \sqrt{f+1}\leq\sqrt{P^{\ast}_{\theta}f+1} \sqrt{P^{\ast}_{\theta}1}, \ {E}_0 \rm{-a.e.}.
\end{equation}
 In contrast with \eqref{E2}, which implies $P^{\ast}_{\theta}1= 1$, ${E}_0$-a.e., then $P^{\ast}_{\theta}(x, \cdot)$ is a probability measure.
Therefore, if and only if $f$ is a constant  under $P^{\ast}_{\theta}(x, \cdot)$  for ${E}_0$-a.e. $x$, the  equation in \eqref{holdereq} holds.

 Moreover, since $p(x,y)>0$ and for some $K>0$, it holds that $\frac{1}{K}p_{\theta_1}(x,y)\leqslant p_{\theta_2}(x,y)\leqslant K p_{\theta_1}(x,y)$, $\theta_1,\theta_2 \in \Theta$, then $p_\theta(x,y)>0$. Therefore, $P^{\ast}_{\theta} 1_A=0$ implies  ${E}_0[1_A]=0$, so
 $f$ is a constant under  ${E}_0$-a.e., which together with $E_0[f]=1$ implies that $f=1,  {E}_0$-a.e..
 \item [(4)] Since $\Phi$ is an increasing function, we have
 \begin{equation}\label{E4}
\Phi(P_\theta f(x))\leq \Phi(\bar{P}f(x))\leq \bar{P}\Phi(f(y)) e^{\Psi(x,y)}
\leq E[p(y,\cdot)\Phi(f)(\cdot)]e^{\Psi(x,y)}.
\end{equation}
From (2),   there exists a $p_\theta$ on $\mathbb{R}^d\times \mathbb{R}^d$, sush that ${P}_\theta f(x)= {E}[p_\theta(x,\cdot) f(\cdot)]$.

  Taking $f=n\wedge \Phi^{-1}\left( {p(x,\cdot)}\right) $ in  \eqref{E4}, we have
\begin{equation*}
e^{-\Psi(x,y)}\Phi(E[p_\theta(x,\cdot)(n\wedge \Phi^{-1}\left( {p(x,\cdot)}\right) )])
\leq  E[p(y,\cdot)\Phi(n\wedge \Phi^{-1}\left( {p(x,\cdot)}\right) )]<\infty.
\end{equation*}
  Letting $n\to\infty$, by monotone convergence theorem, we get
\begin{equation}\label{E5}
 E[p(y,\cdot) p(x,\cdot) ]
 \geq e^{-\Psi(x,y)}\Phi(E[p_\theta(x,\cdot)(\Phi^{-1}\left( {p_\theta(x,\cdot)}\right) )]).
\end{equation}
Since $r\Phi^{-1}(r)$ is convex for $r\geqslant0$, by Jensen's inequality, we have
\begin{equation}\label{E6}
E[p_\theta(x,\cdot)(\Phi^{-1}\left( {p_\theta(x,\cdot)}\right) )]\geq\Phi^{-1}(1).
\end{equation}
Combining \eqref{E5} and \eqref{E6}, we obtain
\begin{equation*}
 E[p(y,\cdot) p(x,\cdot) ]\geq e^{-\Psi(x,y)}.
\end{equation*}
 \item [(5)]  This result is obvious, we  omit it here.

  \end{enumerate}
  Next, we give  an example which  has at most one invariant linear expectation by means  Harnack inequality.
\begin{example}\rm \label{ex3.6}
For $\theta \in \left\{\frac{1}{2},1\right\},$ let $W_t$ be the stand 1-dimensional Brownian motion, consider
 the following SDE,
\begin{equation*}
dX_t=-\theta X_tdt+\sqrt{2\theta}dW_t,\ X_0=x.
\end{equation*}
Then $X_t=e^{-\theta t}x+\int_0^t\sqrt{2\theta}e^{-\theta (t-s)}dW_s$, $X_t \to N(0,1)=\mu_0$ in distribution as $t\to \infty.$ 

Let
\begin{equation*}
P_\theta f=\int_R\frac{1}{\sqrt{2\pi(1-e^{-2\theta })}}\exp\left\{{\frac{{-(z-e^{\theta }x)^2}}{2(1-e^{-2\theta })}}\right\}f(z)dz.
\end{equation*}
Therefore, for every $P_\theta$, there is a same invariant linear expectation $E$, where $E[f]=\int_{\mathbb{R}}fd\mu_0.$
Thus,
\begin{equation*}
E[\bar{P}f]=E\left[\sup_{\theta=\frac{1}{2},1} P_\theta f\right]\leq E\left[P_{\frac{1}{2} }f+P_1f\right]=2E[f],
 \ 0\leq f \in \mathcal{B}_b(\mathbb{R}).
\end{equation*}
Then $E$ is  a  quasi-invariant linear expectation of $\bar{P}$.

Moreover, for every $P_\theta$, $\theta = \frac{1}{2},1$ and $\alpha >1$, it holds that
\begin{equation*}
({P_\theta}f(x))^\alpha \leq {P_\theta}(f(y))^\alpha \exp\{C(\alpha,\theta)|x-y|^2\}\leq \bar{P}(f(y))^\alpha  \exp\{C(\alpha,\theta)|x-y|^2\}.
\end{equation*}
 Then,
 \begin{equation*}
(\bar{P}f(x))^\alpha  \leq \bar{P}(f(y))^\alpha  \exp \left\{\left(C\left(\alpha,\frac{1}{2}\right)+C(\alpha,1)\right)|x-y|^2 \right\},  \ 0\leq f \in \mathcal{B}_b(\mathbb{R}).
\end{equation*} By Theorem \ref{proposition3},  $\bar{P}$ has at most one  invariant linear  expectation.

\end{example}

\subsubsection{Application of Shift Harnack Inequalities}
   Now, we consider the following shift Harnack inequality
    \begin{equation} \label {shiftap1}
\Phi(\bar{P}f(x)) \leq\bar{P}\{\Phi \circ f(e+\cdot)\}(x) e^{C_\Phi(x,e)},  \  0<f \in \mathcal{B}_{b}(\mathbb{R}^d),
\end{equation}
for some $x,e \in \mathbb{R}^d$, and  constant $e^{C_\Phi(x,e)}\geq0$.

\begin{theorem}\label{thm4.2.1} Let $\Phi:[0,\infty)\to [0,\infty)$ be a strictly increasing continuous function such that
\eqref{shiftap1} holds.
 If there exists  a countable set $\Theta_0$ and a family of positive function $\{g_\theta\}_{\theta \in \Theta_0}$ on $\mathbb{R}^d$, such that
$$ \bar{P}f(x)\leq \sup_{\theta\in\Theta_{0}}P_{\theta}(g_\theta f)(x), \ 0<f \in \mathcal{B}_{b}(\mathbb{R}^d), \ x\in \mathbb{R}^d.$$
Then
$\bar{P}$ has a $\sup$ transition density $p(x, y)$ with respect to the Lebesgue
measure.

\end{theorem}
\textbf{Proof.}
 From \eqref{shiftap1}, we have
  \begin{equation}\label{shap2}
\Phi(\bar{P}f(x)) e^{-C_\Phi(x,e)} \leq\bar{P}\{\Phi \circ f(e+\cdot)\}(x).
\end{equation}
Integrating both sides in  \eqref{shap2} with respect to $de$,  it holds that
 \begin{eqnarray}\label{shap3}
\Phi(\bar{P}f(x)) \int_ {\mathbb{R}^d}e^{-C_\Phi(x,e)}de &\leq& \int_ {\mathbb{R}^d}\bar{P}\{\Phi \circ f(e+\cdot)\}(x)de.
\end{eqnarray}
 Moreover, for any Lebesgue-null set $A$, it holds that
 $$\int_ {\mathbb{R}^d}{P}_\theta\{g_\theta \Phi \circ 1_A(e+\cdot)\}(x)de=0.$$
 Then
  \begin{equation*}
  {P}_\theta\{g_\theta\Phi \circ 1_A(e+\cdot)\}(x)=0, \ \rm{ for \ a.e.} \ e.
   \end{equation*}
  Since $\Theta_0$ is a countable set,  thus $\bar{P}\{\Phi \circ 1_A(e+\cdot)\}(x)=0,$  for  a.e. $e$.

Therefore,
  \begin{equation} \label{shiftap4}
  \int_ {\mathbb{R}^d}\bar{P}\{\Phi \circ 1_A(e+\cdot)\}(x)de=0.
  \end{equation}
 By the strictly increasing  properties, we get $\Phi^{-1}(0)=0.$
 Applying $f=1_A$ in \eqref{shap3}, combining with \eqref{shiftap4}, we get
 \begin{equation}
 \bar{P}1_A(x)\leq \Phi^{-1}\left(\frac{ \int_ {\mathbb{R}^d}\bar{P}\{\Phi \circ 1_A(e+\cdot)\}(x)de}
{\int_{\mathbb{R}^d}e^{-C_\Phi(x,e)}de}\right)=0.
\end{equation}
Then for any Lebesgue-null set $A$, we have ${P_\theta}1_A(x)=0$,
 which implies that there exists a density function $p_\theta$ on $\mathbb{R}^d\times \mathbb{R}^d$, such  that
 $${P_\theta}f(x)=\int_{\mathbb{R}^d}p_\theta(x, y)f(y)dy.$$
Thus,  $$\bar{P}f(x)= \sup_{\theta \in \Theta}P_\theta f(x)\leq \int_{\mathbb{R}^d}p(x,y)f(y)dy,$$
 where the $\sup$ density function $p(x, y)=\sup_{\theta \in \Theta}p_\theta(x, y)$.

\begin{example}\rm
Consider the Example \ref{ex3.6} again.
For every $P_\theta$, $\theta \in \{\frac{1}{2},1\}$ and $\alpha >1$, it holds that
\begin{equation*}
({P_\theta}f(x))^\alpha \leq {P_\theta}(f^\alpha(v+\cdot)) (x)\exp\{C(\alpha,\theta)|x-y|^2\}\leq \bar{P}(f^\alpha(v+\cdot)) \exp\{C(\alpha,\theta)|x-y|^2\}.
\end{equation*}
 Thus, there holds the shift Harnack inequality
 \begin{equation*}
(\bar{P}f(x))^\alpha  \leq \bar{P}(f^\alpha(v+\cdot))  \exp \left\{\left(C\left(\alpha,\frac{1}{2}\right)+C(\alpha,1)\right)|x-y|^2 \right\},  \  f \in \mathcal{B}_b^+(\mathbb{R}).
\end{equation*}
 Let $\bar{P}f=\sup_{\theta=\frac{1}{2}, 1}P_\theta f, f\in \mathcal{B}_b^+(\mathbb{R}),$ where
\begin{equation*}
P_\theta f(x)=\int_R\frac{1}{\sqrt{2\pi(1-e^{-2\theta })}}\exp\left\{{\frac{{-(z-e^{\theta }x)^2}}{2(1-e^{-2\theta })}}\right\}f(z)dz, \ \theta=\frac{1}{2}, 1.
\end{equation*}
Since
\begin{eqnarray*}
&&p_{\frac{1}{2}}(x, y)+p_1(x,y)\\
&=&\frac{1}{\sqrt{2\pi(1-e^{-1})}}\exp\left\{{\frac{{-(y-e^{\frac{1}{2} }x)^2}}{2(1-e^{-1 })}}\right\}+\frac{1}{\sqrt{2\pi(1-e^{-2})}}\exp\left\{{\frac{{-(y-ex)^2}}{2(1-e^{-2 })}}\right\}\\
&\leq& \frac{1}{\sqrt{2\pi(1-e^{-1})}}\exp\left\{{\frac{{-(y-e^{\frac{1}{2} }x)^2}}{2(1-e^{-1 })}}+{\frac{{-(y-ex)^2}}{2(1-e^{-2 })}}\right\}.
\end{eqnarray*}
Let $$p(x,y)=\frac{1}{\sqrt{2\pi(1-e^{-1})}}\exp\left\{{\frac{{-(y-e^{\frac{1}{2} }x)^2}}{2(1-e^{-1 })}}+{\frac{{-(y-ex)^2}}{2(1-e^{-2 })}}\right\}.$$
From Theorem \ref{thm4.2.1}, we know that the
 $p(x, y)$ is  a $\sup$ transition density of $\bar{P}$  with respect to the Lebesgue measure.

\end{example}

\paragraph{Acknowledgement.} The author would like to thank Professor Feng-Yu Wang for  guidance and helpful comments, as well as Xing Huang for corrections.

\beg{thebibliography}{99}

 \bibitem{Aida}  S. Aida, Uniformly positivity improving property, Sobolev inequalities and spectral gaps, J. Funct. Anal. 158 (1998) 152--185.

  \bibitem{AK}  S. Aida,    H. Kawabi, Short time asymptotics of certain infinite dimensional diffusion process, Stochastic Analysis and Related Topics, VII, Kusadasi, Japan  (1998) 77--124.

  \bibitem{AZ} S. Aida,   T. Zhang, On the small time asymptotics of diffusion processes on path groups, Potential Anal. 16 (2002) 67--78.

 \bibitem{Bobkov}  S. Bobkov,  I. Gentiland M. Ledoux , Hypercontractivity of Hamilton-Jacobi equations,
J. Math. Pures Appl. 80 (2001)  669--696.

\bibitem{15} L. Denis, M. Hu, S. Peng, Function spaces and capacity related to a sublinear expectation: application to $G$-Brownian motion pathes, Potential Anal. 34 (2011) 139--161.

\bibitem{Gong}  F.  Gong,  F.-Y. Wang,   Heat kernel estimates with application to compactness of manifolds, Q. J. Math. 52 (2001) 171--180.

 \bibitem{Hu} M. Hu, S. Ji, S. Peng, Y. Song, Comparison theorem, Feynman-Kac formula and Girsanov transformation for BSDEs driven by $G$-Brownian motion, Stochastic Process. Appl. 124 (2014), 1170--1195.

\bibitem{HWZ} M. Hu, H. Li,   F. Wang, G. Zheng,  Invariant and ergodic nonlinear expectations for $G$-diffusion processes, Electron. Commun. Probab. 20 (2015) 15 pp.

\bibitem{Kawabi}  H. Kawabi, The parabolic Harnack inequality for the time dependent Ginzburg-Landautype
SPDE and its application, Potential Anal.   22  (2005) 61--84.

\bibitem{Li} X. Li, S. Peng, Stopping times and related It\^{o}'s calculus with $G$-Brownian motion, Stochastic Process, Appl. 121 (2011)  1492--1508

\bibitem{Osuka} E. Osuka, Girsanov's formula for $G$-Brownian motion, Stochastic Process. Appl. 123 (2013) 1301--1318.

\bibitem{peng2} S. Peng, $G$-Brownian motion and dynamic risk measures under volatility
uncertainty, (2007) arXiv: 0711.2834v1.

\bibitem{peng1} S. Peng, $G$-expectation, $G$-Brownian motion and related stochastic
calculus of It\^{o} type, in: Stochastic Analysis and Applications, in: Abel Symp., vol. 2, Springer, Berlin, 2007, pp.541--567.

\bibitem{peng4} S. Peng,  Nonlinear expectations and stochastic calculus under uncertainty-with
robust central limit theorem and $G$-Brownian motion, (2010) arXiv:1002.4546v1.

\bibitem{RW1} M.  R\"ockner and F.-Y. Wang, Supercontractivity and ultracontractivity for (nonsymmetric)
diffusion semigroups on manifolds, Forum Math. 15 (2003)  893--921.

\bibitem{RW2} M.  R\"ockner,  F.-Y. Wang, Harnack and functional inequalities for generalized Mehler
semigroups, J. Funct. Anal.   203 (2003) 237--261.

\bibitem{W7} F.-Y. Wang, Functional inequalities, semigroup properties and spectrum estimates, Infin. Dimens. Anal. Quantum Probab. Relat. Top. 3 (2000) 263--295.

\bibitem{W1} F.-Y. Wang, Harnack inequalities for stochastic partial differential equations,
Springer Briefs in Mathematics, Springer, New York, 2013, pp, ISBN: 978--1--4614--7933-8, 978--1--4614--7934--5.

\bibitem{W} F.-Y. Wang, Logarithmic Sobolev inequalities on noncompact Riemannian manifolds, Probab. Theory Related Fields 109 (1997) 417--424.

\bibitem{W5} F.-Y. Wang, Harnack inequalities for log-Sobolev functions and estimates of log-Sobolev
constants, Ann. Probab.  27 (1999)  653--663.

\bibitem{W6} F.-Y. Wang, Logarithmic Sobolev inequalities: conditions and counter examples, J. Operator
Theory 46  (2001) 183--197.

\bibitem{Xu} J. Xu,  H. Shang,  B. Zhang, A Girsanov type theorem under $G$-framework,
Stoch. Anal. Appl. 29 (2011) 386--406.

\end{thebibliography}

\end{document}